\documentstyle[12pt]{article}
\begin{document}
\newtheorem{thm}{Theorem}
\newtheorem{sbsn}{Subsection}[section]
\title
{Regularity of the distance function to
the boundary}
\author{YanYan Li\thanks{Partially
 supported by
  NSF grant DMS-0401118.}
\\
Department of Mathematics\\
Rutgers University\\
110 Frelinghuysen Road\\
Piscataway, NJ 08854\\
\\
Louis Nirenberg\\
Courant Institute\\
251 Mercer Street \\
New York, NY 10012}

\date{}
\maketitle

\centerline{\it Dedicated to  the memory of Luigi Amerio
}

\input { amssym.def}

\begin{abstract}
Let $\Omega$ be a domain in a smooth complete Finsler manifold,
and let $G$ be the largest open subset of $\Omega$
such that for every $x$ in $G$ there is a unique
closest point from $\partial \Omega$ to $x$ (measured in the Finsler metric).
We prove that the distance function from
$\partial \Omega$
is in  $C^{k,\alpha}_{loc}(G\cup \partial \Omega)$,
$k\ge 2$ and $0<\alpha\le 1$,
 if $\partial \Omega$ is in
 $C^{k,\alpha}$.
\end{abstract}

\setcounter {section} {0}

\section{Introduction}

In \cite{LN} we studied the singular
set of viscosity solutions of some
Hamilton-Jacobi equations.  This was
reduced to the study of the singular
set of the distance function to the 
boundary of a domain $\Omega$ --- for a Finsler
metric.  The singular set was defined as the complement of the following
open set
\begin{equation}
\begin{array}{rl}
G :=& \mbox{the largest open subset of}\ \Omega\
\mbox{such that for every}\ x\ \mbox{in}\ G\
\mbox{there is a}\\
&\mbox
{unique closest point from}\ \partial \Omega\ \mbox{to}\ x\
\mbox{(measured in the Finsler metric).}
\end{array}
\label{2-1}
\end{equation}

In \cite{LN} we stated that if $\partial \Omega$
is in $C^{k,\alpha}$, $k\ge 2$
and $0<\alpha\le 1$, then the distance function from
the boundary belongs to 
$C^{k-1, \alpha}(G\cup \partial \Omega)$.  Recently
Joel Spruck asked to see the proof for
a Riemannian metric and pointed out that
the result would imply that the distance function would
be in $C^{k,\alpha}(G\cup \partial \Omega)$.  In this
paper we provide a proof of that in the
Finsler case.
This paper can be regarded as an addendum to \cite{LN}.

We present two proofs of the $C^{k, \alpha}$ result.
We use the notation, as in \cite{LN},
$$
\int_0^T \varphi(\xi(t); \dot \xi(t))dt
$$
for the length of a curve $\xi(t)$.  
$\varphi(\xi; v)$ is homogeneous of degree one in
$v$.  For fixed $\xi$, the level surface $\varphi(\xi; v)=1$ is smooth,
closed, strictly convex, with positive
principal curvatures.

The first proof uses very little of \cite{LN} and is essentially
self-contained.  The second proof uses some structure from \cite{LN}, and
may be of some interest to some readers.

We actually prove a more general result here, involving
conjugate points from the boundary.

\medskip

\noindent{\bf Definition.  Conjugate Point.}\
Consider a point $y$ on $\partial \Omega$, and consider the
geodesic $\xi(y,s)$ from $y$ going inside $\Omega$ ``normal''
to $\partial \Omega$  (explained below)
with $s$ as arclength.  The conjugate point to $y$ is the first point
$\bar x$ on the normal geodesic such that any point $x''$ on the
geodesic beyond $\bar x$ has,
in any neighborhood of the geodesic, a shorter join
from $\partial \Omega$ to $x''$ than our normal geodesic to $x''$.

\medskip

\noindent{ \bf Normal.}\ A geodesic $\Gamma$ from 
a point $y\in \partial\Omega$ is ``normal'' to
$\partial\Omega$ if for $x$ on $\Gamma$ close to $y$, the geodesic is the shortest
join from $\partial \Omega$ to $x$.

To obtain the regularity in $G$ we prove a slightly more general result
which is local on $\partial \Omega$. Namely, suppose $C$ is a neighborhood on $\partial \Omega$
 of a
point $y$  and that the normal geodesic $\Gamma$ from $y$ to a point $X$
in $\Omega$ is the unique shortest join from $C$ to $X$. If the conjugate point to
$y$ is beyond $X$ then there is a neighborhood $A$ of $X$ such that the distance
from $C$  to any point  in $A$ belongs to $C^{k,\alpha}$.
See Theorem \ref{thm1} below.

We shall make use of special coordinates introduced 
in section 3 of \cite{LN} about a given 
normal geodesic $\Gamma$, going
from a point $y\in \partial \Omega$ into $\Omega$.
In these coordinates $y$ is the origin and the $x_n-$axis
is normal to $\partial \Omega$ there and is the geodesic $\Gamma$.  Furthermore, in
these coordinates, $\varphi$ has the following properties, see
(4.1)-(4.6) in \cite{LN}.
Here 
 Greek letters $\alpha, \beta$ range from $1$ to
$n-1$, and Latin letters $i, j$ range from $1$ to $n$.

\begin{eqnarray}
\varphi(te_n; e_n)&\equiv& 1, 
\label{4-1}
\\
\varphi_{\xi^j}(te_n; e_n)&\equiv& 0, \label{4-2}
\\
\varphi_{v^\alpha}(te_n; e_n)&\equiv& 0, \qquad
\varphi_{v^n}(te_n; e_n)\equiv 1, \label{4-3}
\\
\varphi_{\xi^jv^k}(te_n; e_n)&\equiv& 0,  \label{4-4}\\
\varphi_{v^jv^n}(te_n; e_n)&\equiv& 0,  \label{4-5}\\
\varphi_{\xi^j\xi^n}(te_n; e_n)&\equiv& 0.\label{4-6}
\end{eqnarray}

In these coordinates for $y\in \partial\Omega$ near the origin
the geodesic $\xi(y, s)$ from $y$
``normal'' to $\partial \Omega$ there satisfies 
$$
\dot \xi(y, 0)=V(y)
$$
where $V(y)$ is the unique vector-valued function on $\partial\Omega$
satisfying (here $\nu(y)$ is the Euclidean interior unit normal to
$\partial\Omega$ at $y$)
\begin{equation}
\left\{
\begin{array}{l}
V(y)\cdot \nu(y)>0\\
\varphi(y; V(y))=1\\
\nabla_v\varphi(y; V(y))\
\mbox{is parallel to}\ \nu(y).
\end{array}
\right.
\label{1.16}
\end{equation}

Using these special coordinates, near the origin, $\partial\Omega$ has the form
\begin{equation}
y=(x', f(x')),\quad
f(0')=0,\ \ \nabla f(0')=0'.
\label{10-1}
\end{equation}
We assume that $f\in C^{k, \alpha}$, $k\ge 2$,
$0<\alpha\le 1$.

The result we prove is
\begin{thm}
Assume that the conjugate point of the origin on the geodesic 
$\Gamma=\{te_n\}$ is beyond $e_n$, and that there exists
a neighborhood $C$ of $0'$ on
$\partial \Omega$  such
that $\{te_n\ |\  0\le t\le 1\}$
 is the unique shortest geodesic from $C$ to $e_n$.
Then there exist neighborhoods $A$ of $e_n$ and
${\cal A}$ of $0'$ on $\partial \Omega$
such that for any $X$ in $A$ there is a unique $y\in 
{\cal A}$ and geodesic from $y$ to $X$ which is the 
shortest join
from ${\cal A}$ to $A$.  Furthermore,
if $d(X)$ is its length, then the Jacobian of the map
$X\to (d,y)$ is nonsingular at $e_n$, and
$d$ lies in $C_{loc}^{k,\alpha}$ in $A$, and $y$ lies
in $C_{loc}^{k-1, \alpha}$ in $A$.
\label{thm1}
\end{thm}

\section{Second Variation}

Consider one parameter family of curves $\tau(\epsilon, t)$ from
${\cal A}$ to $\bar te_n$, $\bar t>0$, with $\tau(0,t)=te_n$.   We look at the second 
variation of its length $I[\tau(\epsilon, \cdot)]$.
For $\bar t$ small, it is clearly positive definite.  The first
$\hat t$ for which it fails to be strictly positive definite is the conjugate
point.  For if $\tilde t=\hat t+\delta$, $\delta>0$, then 
the second variation of curves to $\tilde t$ cannot be 
semipositive definite, and there would then
be a shorter connection from ${\cal A}$ to
$\tilde t e_n$ near $\Gamma$.

The standard computation of second variation yields
$$
\frac{d^2}{d\epsilon^2}I[\tau(\epsilon, \cdot)]\bigg|_{ \epsilon=0}
=J\left( \tau_\epsilon\big|_{\epsilon=0} \right)
-f_{x_\alpha x_\beta}(0')\tau^\alpha_\epsilon(0, 0)\tau^\beta_\epsilon(0,0).
$$
Here $J$ is the usual expression
of the second variation if the bottom point were kept at the origin.  Namely,
\begin{equation}
J\left( \tau_\epsilon\big|_{\epsilon=0} \right)
=\int_0^{\bar t}
\left\{ \varphi_{\xi^\alpha \xi^\beta}(te_n; e_n)
\tau^\alpha_\epsilon(0,t)\tau^\beta_\epsilon(0,t)
+\varphi_{v^\alpha v^\beta}(te_n; e_n)
{\dot \tau}^\alpha_\epsilon(0,t) {\dot \tau}^\beta_\epsilon(0,t)
\right\}dt.
\label{13-1}
\end{equation}
Note that $\tau^n_\epsilon$ and ${\dot \tau}^n_\epsilon$ do not occur in $J$.

\section{First proof of Theorem \ref{thm1}}

\subsection{}  Recalling
(\ref{10-1}) we shall denote the normal geodesic from
$y=(x', f(x'))$ by $X=\xi(x', s)$; this is a slight
change of notation.  The geodesic
$\xi$ and $\xi_s$ depend smoothly
on $s$ and their initial data, while the initial data depend
$C^{k-1, \alpha}$ on $x'$.
To prove the theorem, it suffices to show that the 
Jacobian of the mapping $(x',s)$
to $X$ at $(0', 1)$ is nonsingular.
It follows that $d$ and $y$ belong to $C^{k-1,\alpha}$.
Since $\nabla_X d=X_s$, it
follows that $\nabla_Xd$ is in $C^{k-1, \alpha}$ and 
hence $d$ is in $C^{k,\alpha}$ --- as Spruck pointed out to us.

We now prove the Jacobian is nonsingular.

Write $X=(X', X^n)$.  Since 
$X_s(0', 1)=(0', 1)$, the
Jacobian of the mapping $(x', s)$ to $X$ at 
$(0', 1)$ is simply 
$$
M:= \frac{\partial X'}{ \partial x'}(0', 1).
$$
Assume $M$ is singular, without loss
of generality we may suppose that
\begin{equation}
X'_{x_1}(0', 1)=0'.
\label{15-1}
\end{equation}
We construct a perturbation
$\tau(\epsilon, t)$ of $\Gamma=
\{te_n\ |\ 0\le t\le 1\}$ such that
$\displaystyle{
\zeta(t):= \tau_\epsilon\big|_{\epsilon=0} }$ satisfies
\begin{equation}
J[\zeta]=f_{x_\alpha x_\beta}(0')\zeta^\alpha(0)\zeta^\beta(0).
\label{16-1}
\end{equation}

\subsection{}  Consider the geodesic $\xi(\delta e_1, t)$
of length $1$ starting at
$(\delta e_1, f(\delta e_1))$, $0<\delta$ small
and ``normal'' to $\partial \Omega$ there.  Set
\begin{equation}
\zeta(t)=\frac{\partial }{\partial \delta}
\xi(\delta e_1, t)\bigg|_{ \delta=0}.
\label{16-2}
\end{equation}
By (\ref{15-1}),
\begin{equation}
\zeta(1)=0.
\label{16-3}
\end{equation}

We obtain an equation for $\zeta(t)$
by differentiating the geodesic equation
$$
\varphi_{\xi^i}=\frac{d}{dt}
\varphi_{v^i}(\xi; \dot \xi)
$$
with respect to $\delta$, and setting
$\delta=0$.  We find
$$
\varphi_{\xi^i\xi^j}(te_n; e_n)
\zeta^j =\frac {d}{dt}\left( \varphi_{v^iv^j}(te_n;e_n)
{\dot \zeta}^j\right).
$$
Here we have used property (\ref{4-4})
of our special coordinates.  By (\ref{4-6}) and
(\ref{4-5}),
\begin{equation}
\varphi_{\xi^\alpha\xi^\beta}\zeta^\beta= \frac {d}{dt}\left( \varphi_{v^\alpha
v^\beta}
{\dot \zeta}^\beta\right).
\label{17-1}
\end{equation}
We have
\begin{equation}
\zeta^\alpha(0)=\delta^\alpha_1.
\label{17-2}
\end{equation}
In addition,
\begin{equation}
{\dot \zeta}(0)=\frac{\partial }{\partial \delta}
{\dot \xi}(\delta e_1, 0)\bigg|_{\delta=0}
= \frac{\partial }{\partial \delta}
  V(\delta e_1)\bigg|_{\delta=0}
=V_{x_1}(0').
\label{17-3}
\end{equation}

By the last formula in (\ref{1.16}) we have
$$
\nabla_v\varphi\left( (\delta e_1, f(\delta e_1)); V(\delta e_1)\right)
\cdot \left(e_1+f_{x_1}(\delta e_1)e_n \right)=0.
$$
Differentiating in $\delta$ and setting $\delta=0$, we find,
using properties of our special coordinates,
\begin{equation}
\varphi_{v^1v^\beta}(0'; e_n)
V^\beta_{x_1}(0')+f_{x_1 x_1}(0')=0.
\label{18-1}
\end{equation}

Now we introduce the perturbation
$\tau(\epsilon, t)$ as follows
$$
\begin{array}{rll}
\tau^\alpha(\epsilon, t)&=& \epsilon \zeta^\alpha(t),\\
\tau^n(\epsilon, t)&=& te_n
+(1-t) f(\epsilon e_1).
\end{array}
$$
The definition of $\tau^n$ is just to ensure that
$\tau(\epsilon, 0)$ lies on $\partial\Omega$.

According to (\ref{13-1}),
$$
J[\tau_\epsilon\big|_{\epsilon=0}]
=\int_0^1 \left\{
\varphi_{\xi^\alpha\xi^\beta}(te_n; e_n)\zeta^\alpha \zeta^\beta+
\varphi_{v^\alpha v^\beta}(te_n; e_n)
{\dot \zeta}^\alpha {\dot \zeta}^\beta\right\}dt.
$$
Integrating the last expression
by parts we find,
using (\ref{17-1}) , (\ref{17-2}),
(\ref{17-3}) and (\ref{18-1}),
\begin{eqnarray*}
J[\tau_\epsilon\big|_{\epsilon=0}]
&=&
\int_0^1 \left\{
\varphi_{\xi^\alpha\xi^\beta}\zeta^\alpha \zeta^\beta
-\zeta^\alpha \frac{d}{dt}\left(\varphi_{v^\alpha v^\beta}
{\dot \zeta}^\beta\right) \right\}
-  \zeta^\alpha(0)
\varphi_{v^\alpha v^\beta}(0'; e_n){\dot \zeta}^\beta(0)
\\
&=& -\varphi_{v^1v^\beta}(0'; e_n)V_{x_1}^\beta(0')
=f_{x_1 x_1}(0')=
f_{x_\alpha x_\beta}(0') \tau^\alpha_\epsilon(0,0)
\tau^\beta_\epsilon(0,0).
\end{eqnarray*}
It follows from Section 2 that the 
second variation is zero.

\vskip 5pt
\hfill $\Box$
\vskip 5pt

\section{Second proof of Theorem \ref{thm1}}

\noindent{\bf Second proof  of Theorem \ref{thm1}}.\
For $X$ near $e_n$ and for small $\sigma'=(\sigma_1, \cdots, 
\sigma_{n-1})\in \Bbb R^{n-1}$,  let $\tau=\tau(\sigma', X)$ be defined by
$\varphi(X; (\sigma', \tau))=1$ and
$\tau(0', e_n)=1$.   Since $\varphi_{v^n}(e_n; e_n)=1$,
by the implicit function theorem, $\tau$ exists as a smooth
function of $(\sigma', X)$ near $(0', e_n)$.

Let, as on page 111 of \cite{LN}, $\eta=\eta(\sigma', X, t)$ be the unique
smooth function of,
with $\psi=\varphi^2$,
$$
\psi_{ \xi^i}(\eta; {\dot \eta})=
\frac{d}{dt} \psi_{v^i}(\eta; {\dot \eta}),\qquad
t\le 1,
$$
satisfying 
$$
\eta(\sigma', X,1)=X,\qquad
{\dot \eta}(\sigma', X,1)=(\sigma', \tau(\sigma', X)).
$$
As explained in the last two lines of page 108 in \cite{LN},
$\eta(\sigma', X,t)$ is a geodesic with $t$ the arclength.

In the special coordinates described in Section 1,
$\partial\Omega$ has the form 
(\ref{10-1}) near the origin with $f\in C^{k,\alpha}$,
$k\ge 2$, $0<\alpha\le 1$.  Since
$\{te_n\ |\ 0\le t\le 1\}$ is the unique
shortest geodesic from $C$ to $e_n$, we know that for
$X$ close to $e_n$, there exists $x'$ close to $0'$ such that
the ``normal geodesic'' starting
from $(x', f(x'))$ will reach $X$ as a shortest join from $C$
to $X$.  It follows that for some $\sigma'$ close to $0'$
and $t$ close to $0$, we have
\begin{equation}
\left\{
\begin{array}
{rll}
\eta(\sigma', X,t)-(x', f(x'))&=&0,\\
{\dot \eta}^\mu(\sigma', X,t)-V^\mu(x')&=&0,
\end{array}
\right.
\label{25-1}
\end{equation}
where $V(x'):=V(x', f(x'))$ satisfies
(\ref{1.16}). Note that $1-t$ is the distance from $C$ to $X$.

To prove Theorem \ref{thm1}, 
we only need to show that the left hand side of (\ref{25-1}),
denoted as LHS, has nonsingular Jacobian
$\displaystyle{
\frac{\partial (LHS)}{\partial (t,\sigma', x')}
}$ at $(t,\sigma', x', X)=(0,0',0', e_n)$.  Indeed,
this would allow the use of the implicit 
function theorem to show that for $X$ close to
$e_n$ and in a neighborhood of
$(0, 0', 0')$, there 
exists a unique $C^{k-1, \alpha}$ solution
$(t, \sigma', x')=(t(X), \sigma'(X), x'(X))$
of (\ref{25-1}).  Thus, Theorem \ref{thm1}
follows as explained at the beginning of Section 3.
 
Clearly,
\begin{displaymath}
\frac{\partial (LHS)}{\partial t}
(0,0',0', e_n)=
\left(
\begin{array}{c}
{\dot\eta}(0',e_n,0)\\
({\ddot \eta}^\mu(0',e_n,0))
\end{array}
\right)
=
\left(
\begin{array}{c}
0'\\
1\\
0'
\end{array}
\right)
,
\end{displaymath}
a $(2n-1)\times 1$ column vector,
\begin{equation}
\frac{\partial (LHS)}{\partial \sigma'}
(0,0',0', e_n)=
\left(
\begin{array}{c}
(\eta_{\sigma_\alpha}(0',e_n, 0))\\
({\dot \eta}^\mu_{\sigma_\alpha}(0',e_n, 0))
\end{array}
\right)
=\left(
\begin{array}{c}
(\eta^\mu_{\sigma_\alpha}(0',e_n, 0))\\
(\eta^n_{\sigma_\alpha}(0',e_n, 0))\\
({\dot \eta}^\mu_{\sigma_\alpha}(0',e_n, 0))
\end{array}
\right),
\label{28-1}
\end{equation}
a $(2n-1)\times (n-1)$ matrix,
\begin{equation}
\frac{\partial (LHS)}{\partial x'}
(0,0',0', e_n)=
\left(
\begin{array}{c}
-I\\
0\\
-\nabla V'(0')
\end{array}
\right),
\label{29-1}
\end{equation}
a  $(2n-1)\times (n-1)$ matrix,
where $I$ is the $(n-1)\times (n-1)$ identity matrix
and $\nabla V':=(V^\mu_{x_\beta})$.
Thus
\begin{equation}
det\left(
\frac{\partial (LHS)}{\partial (t,\sigma', x')}
(0,0',0',e_n)\right)
=(-1)^{n-1}
det\left(
\begin{array}{cc}
\frac{\partial \eta'}{\partial\sigma'}(0',e_n,0)& -I\\
\frac{  \partial {\dot\eta}'} {\partial\sigma'}(0',e_n,0)& -
\nabla  V'(0')
\end{array}
\right),
\label{30-1}
\end{equation}
where $\frac{\partial \eta'}{\partial\sigma'} :=
(\eta^\mu_{\sigma_\alpha})$ and
$\frac{  \partial {\dot\eta}'} {\partial\sigma'}:=
({\dot\eta}^\mu_{\sigma_\alpha})$.

By the last line in (\ref{1.16}),
$$
\nabla_v\varphi((x', f(x')); V(x'))[e_\beta
+f_{x_\beta}(x')e_n]=0,
$$
i.e.
$$
\varphi_{v^\beta}((x', f(x')); V(x'))+
\varphi_{v^n}((x', f(x')); V(x'))
f_{x_\beta}(x')=0.
$$
Differentiating in $x_\alpha$ and setting $x'=0'$ we find, using properties of our special coordinates
(\ref{4-3}), (\ref{4-4}) and (\ref{4-5}),
\begin{equation}
D^2_{v'}\varphi(0';e_n)\cdot 
\nabla  V'(0')+D^2f(0')=0,
\label{32-1}
\end{equation}
where $D^2_{v'}\varphi:=(\varphi_{v^\beta v^\mu})$.

We now evaluate 
$\frac{\partial {\dot\eta}'}{ \partial\sigma'}(0',e_n,0)$ and
$\frac{\partial \eta'}{ \partial\sigma'}(0',e_n,0)$.
It is proved in section 4.4 of \cite{LN}
that there exists a $C^{2,1}$ function $\tilde f$ near $0'$ satisfying
$$
\tilde f(0')=0, \qquad \nabla \tilde f(0')=0',
$$
\begin{equation}
\left( D^2\tilde f(0')-D^2f(0')\right)>0,
\label{33-1}
\end{equation}
\begin{equation}
\eta(\sigma', e_n, 0)=(y', \tilde f(y')),
\label{33-2}
\end{equation}
\begin{equation}
{\dot \eta}^\mu(\sigma', e_n, 0)=\tilde V^\mu(y'),
\label{33-3}
\end{equation}
where $\tilde V(y'):=
\tilde V((y', \tilde f(y')))$ 
is determined by (\ref{1.16}) with $f$ replaced by $\tilde f$, and
$y'=y'(\sigma')$ satisfies
\begin{equation}
det\left( \frac{\partial y'}{\partial \sigma'}(0')\right)\ne 0.
\label{34-1}
\end{equation}
Note that (\ref{34-1}) is given by (4.9) in
\cite{LN}, while (\ref{33-1}) follows from
corollary 4.15 in \cite{LN} together with the fact that
$e_n$ is not a conjugate point.

Differentiating (\ref{33-3}) in $\sigma_\alpha$ and
setting $\sigma'=0'$ we find$$
{\dot \eta}^\mu_{\sigma_\alpha}(0',e_n,0)=
\tilde V^\mu_{y_\beta}(0')\frac{\partial y_\beta}{\partial \sigma_\alpha}(0'),
$$
i.e.
\begin{equation}
\frac{  \partial {\dot\eta}' }{  \partial\sigma'}(0',e_n,0)=
\nabla  \tilde V'(0') \frac{\partial y'}{ \partial \sigma'}(0').
\label{36-1}
\end{equation}

Differentiating (\ref{33-2}) in $\sigma_\alpha$ and
setting $\sigma'=0'$ we find
\begin{equation}
\frac{\partial \eta'}{ \partial \sigma'}(0',e_n,0)
=
\frac{\partial y'}{ \partial \sigma' }(0').
\label{36-2}
\end{equation}
Since
\begin{displaymath}
\left(
\begin{array}{cc}
\frac{\partial y'}{ \partial \sigma'}(0')& -I\\
\nabla  \tilde V'(0')
\frac{\partial y'}{ \partial\sigma'}(0')&
-\nabla V'(0')
\end{array}
\right)
=
\left(
\begin{array}{cc}
I& -I\\
\nabla  \tilde V'(0')&
-\nabla  V'(0')
\end{array}
\right)
\left(
\begin{array}{cc}
\frac{\partial y'}{ \partial\sigma'}(0')&
\quad\\
\quad& I
\end{array}
\right),
\end{displaymath}
we have,  by putting
(\ref{36-1}) and (\ref{36-2}) into
(\ref{30-1}),
\begin{equation}
\det\left(
\frac{  \partial (LHS) }{  \partial (t,\sigma', x')}
(0,0',0',e_n)\right)
=(-1)^{n-1}
\det\left( \frac{\partial y'}{ \partial \sigma'}(0')\right)
\det\left( 
\begin{array}{cc}
I&-I\\
\nabla  \tilde V'(0')&
-
\nabla   V'(0')
\end{array}
\right).
\label{37-1}
\end{equation}

The proof of (\ref{32-1}), applied to
$\tilde f$ instead of $f$,  yields
\begin{equation}
D^2_{v'}\varphi(0'; e_n)
\nabla  \tilde V'(0') +
D^2 \tilde f(0')=0.
\label{38-1}
\end{equation}
Thus, by (\ref{32-1}) and (\ref{38-1}),
\begin{eqnarray*}
\left(
\begin{array}{cc}
I&-I\\
-D^2\tilde f(0')& D^2 f(0')
\end{array}
\right)
&=&
\left(
\begin{array}{cc}
I&-I\\
D^2_{v'}\varphi(0';e_n)\nabla \tilde V'(0')
&
-D^2_{v'}\varphi(0';e_n)\nabla  V'(0')
\end{array}
\right)\\
&
=&
\left(
\begin{array}{cc}
I&\quad\\
\quad & D^2_{v'}\varphi(0';e_n)
\end{array}
\right)
\left(
\begin{array}{cc}
I&-I\\
\nabla \tilde V'(0')
&
-\nabla  V'(0')
\end{array}
\right),
\end{eqnarray*}
and therefore
\begin{equation}
\det\left( D^2 f(0')- D^2\tilde f(0')\right)
=\det D^2_{v'}\varphi(0'; e_n)
\det
\left(
\begin{array}{cc}
I&-I\\
\nabla \tilde V'(0')
&
-\nabla  V'(0')
\end{array}
\right).
\label{39-1}
\end{equation}

Since $D^2_{v'}\varphi(0'; e_n)$ is positive definite, 
we deduce from (\ref{37-1}),
(\ref{34-1}) and (\ref{39-1}) that
\begin{displaymath}
\det
\left(
\frac{   \partial (LHS) }{  \partial(t,\sigma', x')}
(0,0',0',e_n)
\right)
\ne 0.
\end{displaymath}

\vskip 5pt
\hfill $\Box$

\end{document}